\documentclass[12pt,oneside]{article}
\usepackage{graphicx}

\newcommand{\Titolo}{\vskip 0.0cm \fontsize{24pt}{24pt} \selectfont}
\newcommand{\Abstract} {
     \vskip 0.3cm \fontsize{12pt}{12pt} \selectfont
     Abstract
     \vskip 0.3cm
     \fontsize{10pt}{10pt}
     \selectfont
     }
\newcommand{\Report}{\vskip 0.4cm \fontsize{11pt}{11pt} \selectfont}
\newcommand{\myline}{ \noindent \rule{13cm}{0.75pt} }
\newcommand{\Autore}{\vskip 0.0cm \fontsize{12pt}{12pt} \selectfont }
\newcommand{\Indirizzo}{\vskip 0.0cm \fontsize{10pt}{10pt} \selectfont }

\begin{document}

\Titolo
\centerline{        Optimization Strategies in Complex Systems  }
\myline
\Autore
{\noindent          L. Bussolari, P. Contucci, C. Giardin\`a}
\Indirizzo
{                   Dipartimento di Matematica,
                    Universit\`a di Bologna,
                    Piazza di Porta S. Donato 5,
                    40127 Bologna,
                    {\em \{bussolar, contucci, giardina\}@dm.unibo.it}}
\vspace{0.1cm}
\Autore
{\noindent          C. Giberti }
\Indirizzo
{                   Dipartimento di Informatica e Comunicazione,
                    Universit\`a dell'Insubria,
                    Via Maz\-zi\-ni 5,
                    21100 Varese
                    {\em claudio.giberti@uninsubria.it}}
\vspace{0.1cm}
\Autore
{\noindent          F. Unguendoli, C. Vernia }
\Indirizzo
{
                    Dipartimento di Matematica Pura ed Applicata,
                    Universit\`a di Modena e Reggio Emilia,
                    Via Campi 213/b,
                    41100 Modena
                    {\em \{unguendoli, vernia\}@unimore.it}}

\myline

\Abstract
{\it
We consider a class of combinatorial optimization problems that emerge
in a variety of domains among which: condensed matter physics, theory of
financial risks, error correcting codes in information transmissions,
molecular and protein conformation, image restoration. We show the
performances of two algorithms, the``greedy'' (quick decrease along the
gradient) and the``reluctant'' (slow decrease close to the level curves)
as well as those of a``stochastic convex interpolation''of the two.
Concepts like the average relaxation time and the wideness of the
attraction basin are analyzed and their system size dependence
illustrated.
}

\Report

A decision-making problem is often formulated as minimization
of a function of several variables (the {\em cost
function}) possibly subjected
to some constraints. Consider, for example, the following
combinatorial optimization problem
\footnote{from the {\it Millennium Prize Problems} list at
{\it Clay Mathematics Institute}, http://www.claymath.org }:
``Suppose that you are organizing housing accommodations for a group of four
hundred university students. Space is limited and only one hundred of the
students will receive places in the dormitory. To complicate matters, the
Dean has provided you with a list of pairs of incompatible students, and
requested that no pair from this list appear in your final choice''.
It is clear that one may check ``easely" (in a polynomial time) if
a proposed solution is correct; nevertheless the extensive
search of the solution is extremely {\em hard} because
the total number of possible
accommodations is as large as ${400 \choose 100}$,
a number which is larger of the estimated
number of atoms in the universe!
This is an example of a so-called NP complete problem,
i.e. a problem for which no polynomial algorithm is known
to find the solution.
There are many instances of NP complete problems coming for very
different subjects. In the classical formulation of the theory of
financial risk with short-selling included, the choice of the
optimal portfolio from the historical prices correlations is
casted in the selection of an investment strategy among a huge
number of possibilities \cite{BP}. Another example comes from the
theory of error-correcting codes, where one faces the problem of
reconstructing an original sequence of bits which has been
corrupted during transmission \cite{Ni}. More familiar to
physicists is the NP complete problem coming from the realm of
condensed matter. Here particles carrying a magnetic moment (spin)
is modeled by dichotomic variables $\sigma_i = \pm 1$. The
interaction between different atoms is described by a function
$H(J,\sigma)$ (the Hamiltonian or energy), which plays the role of
the cost function, and which depends also on some randomness
through the set of variables $J$. The optimization problem amounts
to find, for a given instance of the problem specified by the
variables $J$, the ground state configuration, i.e. the spin
configuration which minimize energy.

It is our aim in this report to show how new concepts and
techniques from Statistical Mechanics can be helpful in finding
approximated optimal solutions to complex NP complete problems
\cite{MPV}.

We perform a statistical analysis of energy-decreasing
algorithms on a specific mean-field spin model with
complex energy landscape. We specifically address the following question:
in the search of low energy configurations is it convenient
(and in which sense) a quick decrease along the gradient
(greedy dynamics) or a slow decrease close to the level curves
(reluctant dynamics)?
Average time and wideness of the attraction basins are introduced
for each algorithm together with a convex interpolation among the two
and experimental results are presented for different system sizes.
We found that while the reluctant algorithm performs better for a fixed
number of trials, the two algorithms become basically equivalent for a given
elapsed  time due to the fact that the greedy has a shorter relaxation
time which scales linearly with the system size compared to a quadratic
dependence for the reluctant. A final test is also performed in a stochastic
convex combination of the two algorithms: at each step the
motion is greedy with probability $P$ and reluctant
with probability $1-P$. It is found that for large $N$
and for fixed running times a substantial improvement
is obtained with a $P=0.1$

\vspace{0.3cm}
{\bf Model and Algorithms }
\vspace{0.3cm}

\noindent
We consider the paradigmatic model
of complex spin systems, the so-called
Sherrington-Kirkpatrick model \cite{SK}.
This is defined by the following
Hamiltonian
\begin{equation}
H(J,\sigma) = - \frac{1}{2}\sum_{i,j=1}^N J_{ij}\sigma_i\sigma_j
\end{equation}
where $\sigma_i=\pm 1$ for $i=1,\dots,N$ are Ising spin variables
and $J_{ij}$ is an $N\times N$ symmetric matrix which
specifies local interaction between them.
The $J_{ij}$ are independent identically distributed
symmetric gaussian random variables
$(J_{ij}=J_{ji},\;J_{ii} = 0)$ with
zero mean and variance $1/N$,
in order to have a sensible thermodynamic
limit.
Since this is a
disordered model one is interested in the
quenched average ground state energy.
For each $N$ this is defined as:
\begin{equation}
e^{GS}_N = Av \left\{ \frac{1}{N}{\inf}_{\sigma}
H_N(J,\sigma)\right\}
\end{equation}
where we denoted by $Av\{\cdot\}$ the average over
the couplings. Analytical knowledge of this
quantity is available in the thermodynamical limit
$N\rightarrow\infty$ using Parisi Ansatz for
replica symmetry breaking theory:
$e^{GS}_{\infty} = -0.7633$ \cite{MPV},
a result which has been confirmed
by numerical simulations using
finite size scaling, yielding
$e^{GS}_{\infty} = -0.76 \pm .01$
\cite{KS}.

The {\em greedy} and {\em reluctant} dynamics work as follows. The
initial spin configuration at time $t = 0$ is chosen at random
with uniform probability. Then, at each step $t$, the whole
spectrum of energy change obtained by flipping one of the spin is
calculated
\begin{equation}
\Delta E_i = \sigma_i(t)\sum_{j \neq i} J_{ij} \sigma_j(t)
\end{equation}
The configuration is then updated at time $t+1$ by flipping the
spin which corresponds to the largest (greedy) or smallest
(reluctant) energy decrease. Both the dynamics follow an energy
descent trajectory till they arrive to a $1$-spin-flip stable
configuration, i.e. a configuration whose energy can not be
decreased by a single spin-flip. These represent local minima in
energy landscape at zero temperature and the ground state is one
of them. Moreover, we investigate the efficiency of a stochastic
convex combination of the two algorithms: with probability $0\leq
P \leq 1$ we perform a greedy move and with probability $1-P$ the
corresponding reluctant move. The deterministic dynamics are
obtained at $P=1$ (greedy) and $P=0$ (reluctant), respectively.
Intermediate values of $P$ are stochastic dynamics where the
greedy and reluctant moves are weighted by the probability $P$.

\vspace{0.3cm}
{\bf Computational Resources and Results }
\vspace{0.3cm}

\noindent Our simulations used about 5000 hours of CPU time on the
machine IBM SP3. Parallelization has been highly efficient due to
the fact that we run different disorder realization and/or
different initial conditions on each processor, simply averaging
the results at the end of the elaboration.

\begin{figure}[h]
\begin{center}
\includegraphics[width=10.cm,angle=-90]{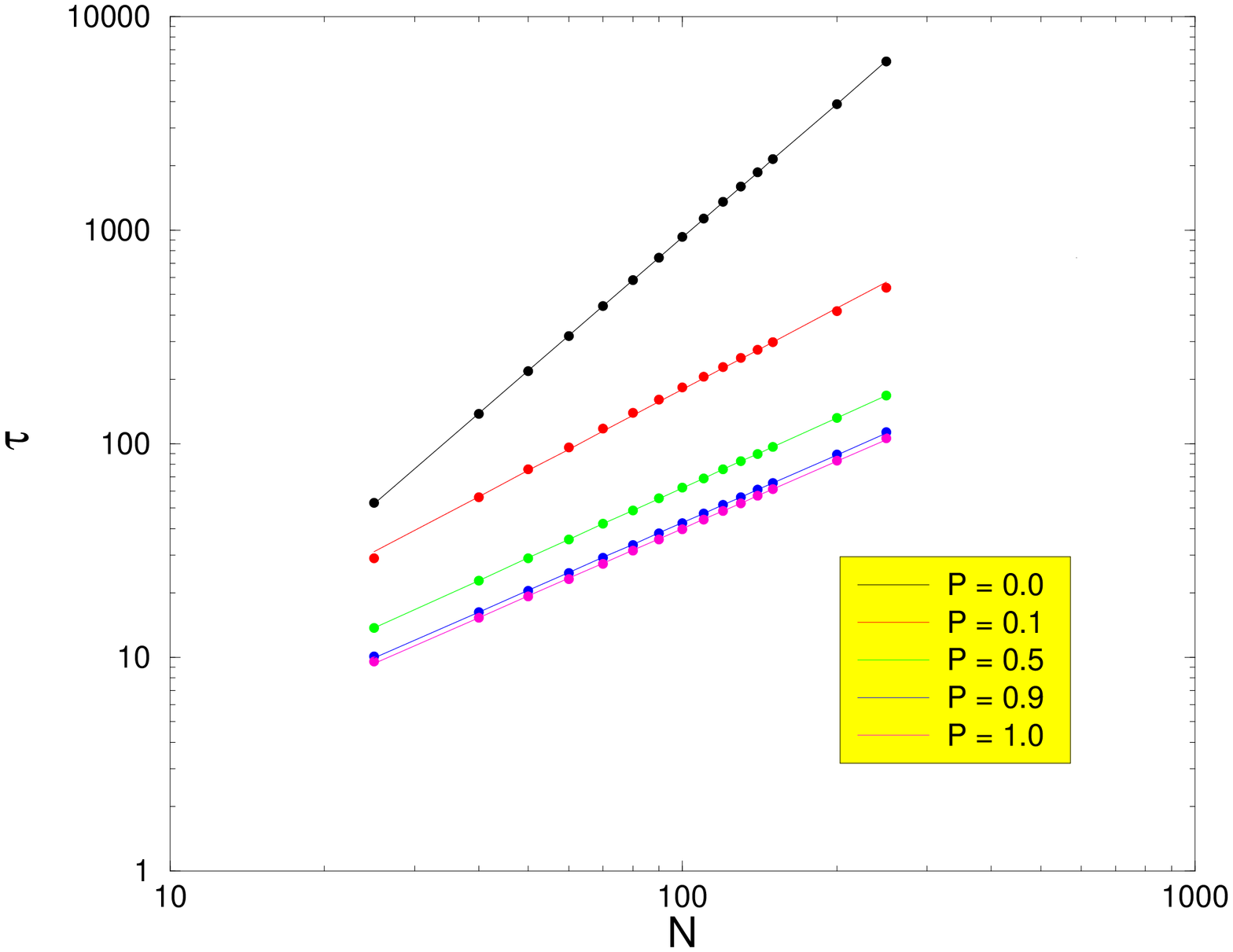}
\caption{
{\small The average time to reach a
metastable configuration for different values
of $P$. Top to bottom: $P = 0$ (reluctant), $P=0.1$, $P=0.5$,
$P=0.9$, $P=1$ (greedy). The continuous lines are the numerical
fits to power law: $\tau(N) \sim N^{\alpha}$,
with $\alpha = 2.07, 1.26, 1.08, 1.05, 1.04$
from top to bottom.}
}
\label{sk-intrinsic}
\end{center}
\end{figure}

First of all, we analyzed the average time $\tau$ of the dynamics
for different values of $P$, which is easily accessible to
measurements and has good self-averaging properties. This is
defined as
\begin{equation}
\label{time}
\tau = \frac{1}{M}\sum_{i=1}^{M} t_i
\end{equation}
where $t_i$, $i = 1,\ldots, M$ is the time of each realization of
the dynamics, measured by counting the number of ``spin flip''
necessary to reach a metastable configuration.. The number of
trials $M$ is an increasing function of the system size. For the
largest sizes $(N = 250,300)$ we used up to $10^9$ initial
configurations. Results are shown in Fig. (\ref{sk-intrinsic}),
together with the best numerical fits. Note the progressive
increase of the slope in log-log scale from an almost linear law
for greedy (bottom) $\tau (N)\sim N^{1.04}$ to an almost quadratic
law for reluctant (top) $\tau (N)\sim N^{2.07}$. However, an
interesting result is that for $P=0.1$ we have still have $\tau
(N)\sim N^{1.26}$, i.e. a stochastic algorithm which makes on
average one greedy move (and nine reluctant moves) out of ten has
a much smaller average time than the deterministic reluctant
algorithm $P=0$. We notice that the exponents for greedy (resp.
reluctant) algorithm are very close to the integers 1 (resp. 2)
with an observed slow crossover between the two for intermediate
$P$. It would be interesting to have a theoretical understanding
of this phenomenon even if only at a heuristic level. We plan to
return over this problem in a future work.

\begin{figure}[h]
\begin{center}
\includegraphics[width=10.cm,angle=-90]{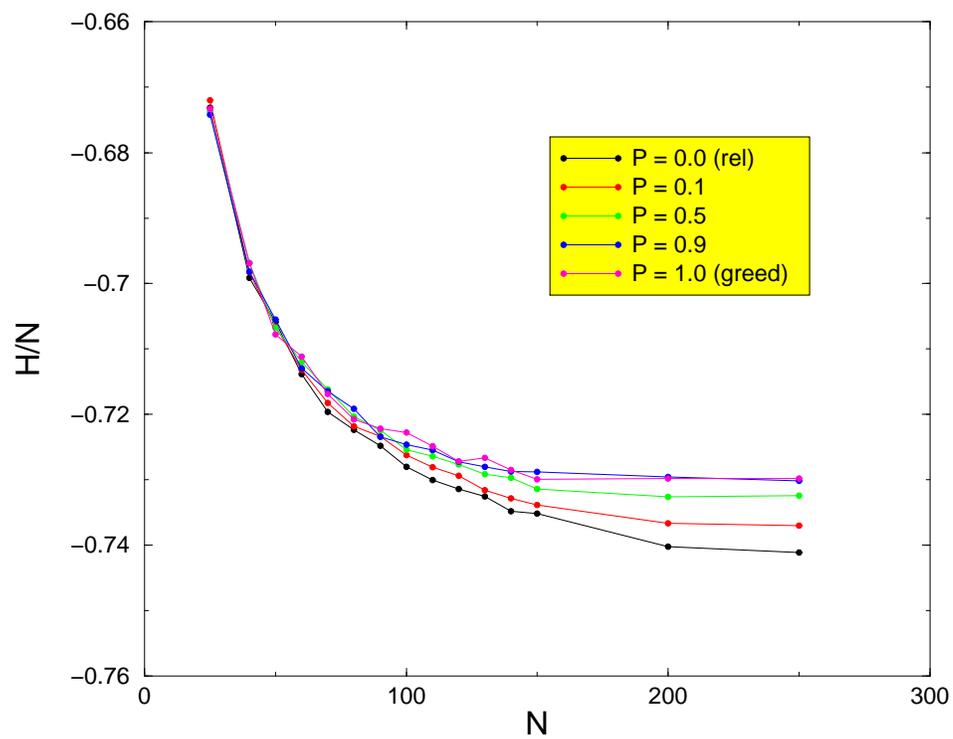}
\caption{
{\small Lowest energy value for
a fixed number of $N$ initial conditions
for different value of $P$. Bottom to top:
$P = 0$ (reluctant), $P=0.1$, $P=0.5$, $P=0.9$, $P=1$ (greedy) }
}
\label{sk-fix-init-cond}
\end{center}
\end{figure}

Next, we measured the lowest energy value found for a fixed number
of initial conditions for different probabilities $P$. One has to
choose a protocol to fix the number of initial conditions.
Obviously, the larger the system size the bigger must be the
number of trials. We tried different choices obtaining similar
results. For the sake of space we show in Fig.
(\ref{sk-fix-init-cond}) the results of the run where we choose
$N$ initial conditions for a system of size $N$. The data have
been averaged on $1000$ disorder realizations. We see that the
smaller is the probability of making greedy moves, the lower is
the energy found. The best result is obtained for $P=0$, which
corresponds to deterministic reluctant dynamics. This means that,
ignoring the total amount of time and imposing constraint only on
the number of initial conditions, reluctant dynamic is the most
efficient in reaching low energy states, i.e. it has a larger
basin of attraction.

\begin{figure}[h]
\begin{center}
\includegraphics[width=10.cm,angle=-90]{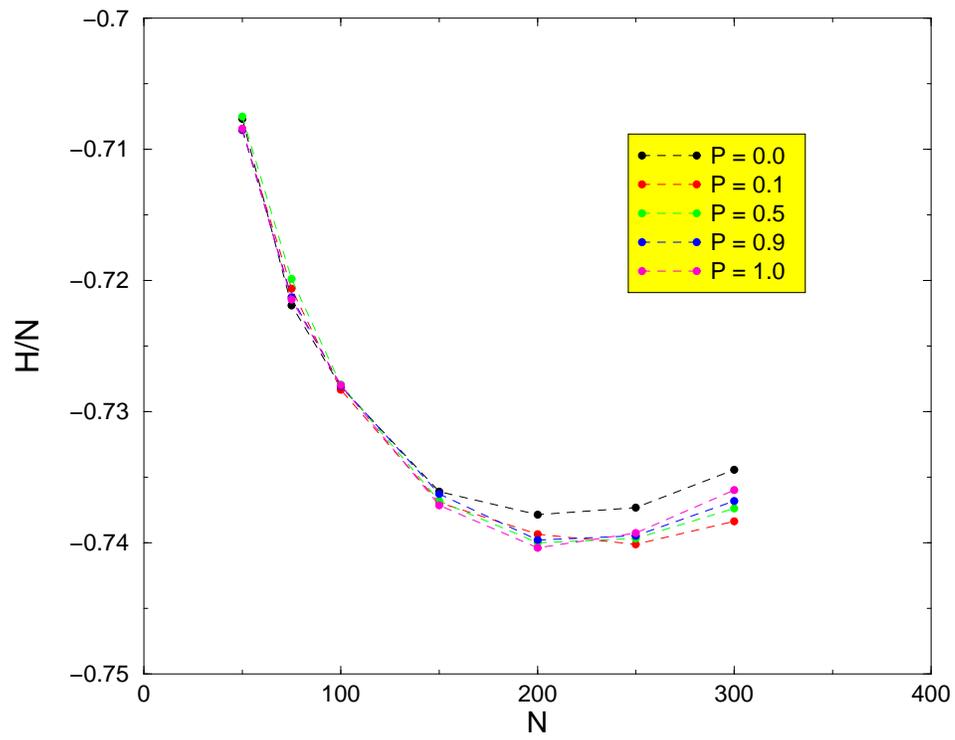}
\caption{
{\small Lowest energy value for a
fixed elapsed time of $100$ hours (each run) on a IBM SP3 for
different value of $P$ (see legend)}
}
\label{sk-fix-elapsed}
\end{center}
\end{figure}

Finally, we compared results of different probabilities in the
case one considers a fixed elapsed time. As an example, we present
results for an elapsed time of $100$ hours of CPU on a IBM SP3 for
$N$ in the range $[50,300]$. We considered again $1000$ disorder
realizations and assigned the same time length to each sample ($6$
minutes). Obviously, in this way reluctant dynamics starts from a
smaller number of initial conditions than greedy, because its
relaxation time is longer. In Fig. (\ref{sk-fix-elapsed}) we
plot the values of the lowest energy state as a function of $N$.
We can see from the data that, for a fixed elapsed time, greedy
dynamics ($P=1$) find lower energy states than reluctant ($P=0$).
Moreover, we observe that the best result is obtained for $P=0.1$.
Thus, we suggest that the more power full strategy to find low
energy state using greedy and reluctant dynamic is a combination
of them, where most of the steps the move is reluctant and on a
small fraction of steps (say $0.1$) the move is greedy.

The results of the present work is extensively presented in ref
\cite{Noi1}. Improvements of the greedy and reluctant algorithms
is presently under study (on IBM SP4), by permitting also {\em
increase} in energy with exponential decrease in time, in the same
spirit of the well-known Simulated Annealing strategies
\cite{Noi2}.

\vspace{0.3cm} {\bf Acknowledgments} We thank the Cineca staff for
technical support and in particular Dr. Erbacci and Dr. Calonaci
for kind assistance. We thank Prof. S. Graffi and Prof. I.
Galligani for all the encouragement and the support of Cineca
grant. \vspace{0.3cm}

\myline


\newpage


\end{document}